\documentclass[preprint,a4paper]{article}

\usepackage{amssymb}
\usepackage{amsfonts}
\usepackage{latexsym}

\newtheorem {theor} {\bf Theorem}
\newtheorem {lemma} {\bf Lemma}

\newtheorem{crl} {\bf Corollary}
\newtheorem{prop} {\bf Proposition}
\newtheorem{dfn}{\bf Definition}
\newtheorem {algr} {\bf Algorithm}

\title {Mathematical Modeling in the Textile Industry}
\author {Krasimir Yordzhev\and Hristina Kostadinova}
\date {}

\begin {document}

\maketitle

\begin{center}
krasimir.yordzhev@gmail.com \hspace{1cm} hkostadinova@gmail.com
\end{center}

\begin{abstract}
A mathematical model, describing some different weaving structures, is made in this article.  The terms self-mirror and rotation-stable weaving structure are initiated here. There are used the properties and operations in the set of the binary matrices and an equivalence relation in this set. Some combinatorial problems about finding the cardinal number and the elements of the factor set according to this relation  is discussed. 
  We propose an algorithm, which solves these problems. The presentation of an arbitrary binary matrix using sequence of nonnegative integers is discussed. It is shown that the presentation of binary matrices using ordered n-tuples of natural numbers makes the algorithms faster and saves a lot of memory.
Implementing these ideas a computer program, which receives all of the examined objects, is created.
  In the paper we use object-oriented programming using the syntax and the semantic of C++ programming language. Some advantages in the use of bitwise operations are shown.
  The results we have received are used to describe the topology of the different weaving structures.
\end{abstract}

{{\bf Key words:}\it   binary matrix, weaving structures, permutation matrix, equivalence relation, C++ programming language, bitwise operations, computer algebra}

{\bf 2010 Mathematics Subject Classification:} 05B20, 15B34, 68N15

\section{Introduction.} It is well  known \cite{borzunov,yordzhev,KYKo}, that the interweaving of the fibres in certain weaving structure can be coded using square \textit{binary} (or (0,1), or  \textit{boolean}) matrix, i.e. all elements of this matrix are 0 or 1. The fabric represented by this matrix exists if and only if in each row and in each column of the matrix there is at least one 1 and at least one 0. Two different matrices correspond to one and the same weaving structure if and only if the first matrix is receive from the second  by several cycle moves of the first row or column to the last place.

Let $n$ be a positive integer. Let us denote by ${\cal B}_n$  the set of all $n\times n$ binary matrices, and let us denote by ${\cal Q}_n$  the set of all $n\times n$ binary matrices which have at least one 0 and at least one 1 in every row and every column. It is obvious, that ${\cal Q}_n \subset {\cal B}_n$. About the necessary definitions and denotes  in the theory of matrices we refer to \cite{sachkov} and \cite{tarakanov}. It is not difficult to see, that
\begin{equation} \label{bn}
\left| {\cal B}_n \right| =2^{n^2}
\end{equation}

If $A=(a_{i\, j} ) \in {\cal B}_n$, then $A^T =(a_{j\, i} )$, $1\le i,j \le n$ denotes the {transposed} matrix $A$.

We are interested in the subset ${\cal P}_n \subset {\cal Q}_n$  of all \textit{permutation} matrices, i.e. binary matrices which have a single 1 in every row and every column. It is well known \cite{tarakanov} that the set ${\cal P}_n$ together with the operation multiplication of matrices is a group, isomorphic to the symmetric group ${\cal S}_n$, where the set
\begin{equation}\label{S_n}
{\cal S}_n =\left\{ \left. \left(
\begin{array}{cccc}
1    &    2    &  \cdots &    n\\
i_1  &    i_2  &  \cdots &    i_n
\end{array}
\right)  \right|  0\le i_k \le n, k=1,2,\ldots ,n, i_k \ne i_l \; \mbox{\textrm{for}} \; k\ne l
\right\}
\end{equation}
is made of all one to one maps of the elements of the set $\{ 1,2,\ldots ,n\}$ to itself. If $M =(m_{k\,l})\in {\cal P}_n$ and the corresponding element in this isomorphism is
$\displaystyle \left(
\begin{array}{cccc}
1    &    2    &  \cdots &    n\\
i_1  &    i_2  &  \cdots &    i_n
\end{array}
\right)
\in {\cal S}_n$, this means that $m_{k\, i_k} =1$, $k=1,2,\ldots,n$ and all other elements of M are equal to 0.

Let $t\in \{1,2,\ldots ,n\}$ and let
$\displaystyle \rho =
\left(
\begin{array}{cccccc}
1    &    2    &  \cdots &  t    & \cdots  &  n  \\
i_1  &    i_2  &  \cdots &  i_t  & \cdots  &  i_n
\end{array}
\right)
\in {\cal S}_n$. We let $(t)\rho =i_t$  denote the image $i_t$ of the number $t$ using the map $\rho$.  $(t)\rho_1 \rho_2 =((t)\rho_1 )\rho_2$ is satisfied   by definition for arbitrary $\rho_1 ,\rho_2 \in S_n$.

As it is well-known \cite{tarakanov}, if we multiply a $n\times n$ matrix $A$ from the right with a permutation matrix $M\in {\cal P}_n$, then this is the same as the changing places of the columns of $A$. In this case, if the corresponding element of  $M\in {\cal P}_n$ in the above-described isomorphism is
$\displaystyle \left(
\begin{array}{cccc}
1    &    2    &  \cdots &    n\\
i_1  &    i_2  &  \cdots &    i_n
\end{array}
\right)
\in {\cal S}_n$,
then after the multiplication we get a matrix $B=AM$ with $k$th column of $B$ equal to $i_k$th column of $A$, $k=1,2,\ldots ,n$. Analogously when we want to exchange the places of rows we multiply $A$ from the left with $M^T$.

Identity element of the group ${\cal P}_n $ is the identity matrix $E_n$, all elements in the leading diagonal are equal to 1 and all elements are equal to 0 everywhere else. The identity element of the group ${\cal S}_n$ is the element $\displaystyle \left(
\begin{array}{cccc}
1    &    2    &  \cdots &    n\\
1    &    2    &  \cdots &    n
\end{array}
\right)$,

We say, that the binary $n\times n$ matrices $A$ and $B$ are equivalent and we write $A\sim B$, if one matrix is made by the other after several consecutive cycle moves of the first row or column to the last place. In other words, if
$A,B\in {\cal Q}_n$ and $A\sim B$, then with the help of these matrices we code one and the same weaving structure. It is obvious, that the relation in the set ${\cal B}_n$ is a equivalence relation. The equivalence class according to the relation $\sim$ with the matrix $A$ we denote $\overline{A} $, and the sets of equivalence classes in ${\cal B}_n$  and   ${\cal Q}_n$ (factor set) according to $\sim$ with $\overline{{\cal B}_n}$ and $\overline{{\cal Q}_n}$. We consider that $\overline{{\cal B}_n}$  and $\overline{{\cal Q}_n}$ are described if there is a random representative of each equivalence class.

The equivalence classes of ${\cal B}_n$ by the equivalence  relation
 $\sim$ are a particular kind of  {\it double coset} \cite{bogopolski,serre,vos}. They make use of substitutions group theory  and linear representation of finite groups theory \cite{serre}.

The elements of $\overline{{\cal Q}_n}$ we call \textit{interweavings}. In that case the number $n$ is called \textit{repeating} of the interweaving. These terms are borrowed from the textile industry.

It is naturally to arise a lot of combinatorial problems, which take place in practice in the weaving industry, connected with the different subsets of $\overline{{\cal Q}_n}$, i.e. with the different classes of interweavings. Some of these classes we examine in the present work.

\section{Some interesting classes of interweavings.}
It is easy to see, that if $A\in {\cal P}_n$ and $B\sim A$, then $B\in {\cal P}_n$. Interweavings which representatives are elements of the set ${\cal P}_n$ of all permutational matrices are called \textit{primary} interweavings. A formula and an algorithm to calculate the number of all primary interweavings with a arbitrary repetition $n$ are shown in \cite{yordzhev,bansko}.

We examine the matrix
\begin{equation}\label{P}
P =
\left(
\begin{array}{ccccc}
0 & 1 &  0   & \cdots & 0 \\
0 & 0 &  1   & \cdots & 0 \\
\vdots    & \vdots    & \vdots   & \vdots  & \vdots   \\
0 & 0 &  0 & \cdots & 1  \\
1 & 0 &  0 & \cdots & 0  \\
\end{array}
\right)
=\left( p_{i\, j} \right) \in {\cal P}_n ,
\end{equation}
where $p_{1\, 2} = p_{2\, 3} = \cdots =p_{i\, i+1} = \cdots =p_{n-1\, n} = p_{n\, 1} = 1$, and all other elements of $P$ are equal  to 0. At the above described isomorphism of the group of permutation matrices with the symmetric group, the matrix $P$ corresponds to the element
\begin{equation}\label{pi}
 \pi=
\left(
\begin{array}{cccccc}
1 & 2 & 3 & \cdots & n-1  &  n\\
2 & 3 & 4 & \cdots & n    &  1
\end{array}
\right)
\in {\cal S}_n .
\end{equation}

It  is not difficult to see that
$P^t \ne E_n $ when $1\le t<n$ and $P^n =E_n $,
where $E_n$ is identity $n\times n$ matrix.
\begin{equation}\label{ttttt}
P^{k+n} =P^k
\end{equation}
for every natural number $k$.

Let $A\in {\cal B}_n$ and let
$$B=PA$$
and
$$C=A P .$$

It is easy to prove (see \cite{tarakanov}), that the first row of $B$ is equal to the second row of $A$, the second row of $B$ is equal to the third row of $A$ and so on, the last row of $B$ is equal to the first row of $A$, i.e. the matrix $B$ is created by the matrix $A$ by moving the first row to the last place, and the other rows are moved one level upper.

Analogously we convince that $C$ is made from $A$ by moving the last column to the first place, and the other column move one position to the right.

Have in mind what is above-described, it is easy to prove the following:

\begin{lemma}\label{l1}
Let  $A,B\in {\cal B}_n$. Then $A\sim B$,  if and only if, there exist integers $k,l\in [0,n-1]$, so that
$$A=P^k B P^l ,$$
 where $P$ is the matrix given by the formula (\ref{P}).

\hfill $\Box$
\end{lemma}

\begin{crl}\label{c2}
All elements of an arbitrary equivalence class of ${\cal B}_n$ according to the relation $\sim$ can be placed in a rectangular $s\times t$ table, where $s$ and $t$ are divisors of $n$.

\hfill $\Box$
\end{crl}

\begin{crl}\label{c1}
Each equivalence class of ${\cal B}_n$ according to the relation $\sim$ contains no more than $n^2$ elements.

\hfill $\Box$
\end{crl}

Except the well-known from the linear algebra operation matrix transposition, here we take more similar matrix operations.

Let
\begin{equation}\label{A}
A=
\left(
\begin{array}{cccc}
a_{1\, 1} & a_{1\, 2} & \cdots & a_{1\, n} \\
a_{2\, 1} & a_{2\, 2} & \cdots & a_{2\, n} \\
\vdots    & \vdots    &        & \vdots    \\
a_{n\, 1} & a_{n\, 2} & \cdots & a_{n\, n} \\
\end{array}
\right)
\in {\cal B}_n
\end{equation}

If $A$ is a square binary matrix, represented by the formula (\ref{A}), then by definition
\begin{equation}\label{A^S}
A^S =
\left(
\begin{array}{cccc}
a_{1\, n} & a_{1\, n-1} & \cdots & a_{1\, 1} \\
a_{2\, n} & a_{2\, n-1} & \cdots & a_{2\, 1} \\
\vdots    & \vdots    &        & \vdots    \\
a_{n\, n} & a_{n\, n-1} & \cdots & a_{n\, 1} \\
\end{array}
\right)
,
\end{equation}
i. e. $A^S$  is created from $A$ as the last column of $A$ becomes the first, the column before last - second and so on, the first column becomes last. In other words, if $A=(a_{i\, j} )$, then $A^S =(a_{i\, n-j+1} )$, $1\le i,j \le n$.

It is obvious, that
$$ \left( A^S \right)^S =A.$$

We say, that the matrix $A\in{\cal B}_n$ is \textit{ a mirror image} of the matrix $B\in {\cal B}_n$, if $A^S =B$.

It is easy to see, that if the matrix $A$ is a mirror image of the matrix $B$, then $B$ is a mirror image of $A$, i.e. the relation ''mirror image'' is symmetric.

In general $A\ne A^S$. If $A = B^S$ and $B = C^S$, then in general we have
$A=B^S =\left( C^S \right)^S =C\ne C^S$.
Therefore, the relation ''mirror image'' is not reflexive and is not transitive.

We examine the matrix
\begin{equation}\label{S}
S =
\left(
\begin{array}{cccccc}
0 & 0 & \cdots & 0 & 0 & 1 \\
0 & 0 & \cdots & 0 & 1 & 0 \\
0 & 0 & \cdots & 1 & 0 & 0 \\
\vdots    & \vdots    &        & \vdots   & \vdots  & \vdots   \\
1 & 0 & \cdots & 0 & 0 & 0 \\
\end{array}
\right)
=\left( s_{i\, j} \right) \in {\cal P}_n ,
\end{equation}
where for each $i=1,2,\ldots ,n$ $s_{i\, n-i+1} =1$ and $s_{i\, j} =0$ as $j\ne n-i+1$. According to the above-described isomorphism of the group of permutational matrices with symmetric group, the matrix $S$ corresponds with the element
\begin{equation}\label{sigma}
\sigma =
\left(
\begin{array}{ccccc}
1 & 2 & 3 & \cdots & n\\
n & n-1 & n-2 & \cdots & 1
\end{array}
\right)
\in {\cal S}_n .
\end{equation}

Obviously $S$ is a symmetric matrix, i.e. $S^T = S$. We check directly, that $S^2 =E_n $.

Is not difficult to see, that
\begin{equation}\label{8}
A^S = A S
\end{equation}
is satisfied for each $A\in {\cal B}_n$.

\begin{lemma}\label{l2}
If $P$ and $S$ are the matrices given by the formulas (\ref{P}) and (\ref{S}), then
\begin{equation}\label{PlS}
P^l S=SP^{n-l}
\end{equation}
for each $l=0,1,2,\ldots ,n-1$.

\end{lemma}

Proof.
Let us denote by $\oplus$ and $\ominus$ the operations  addition and subtraction in the ring ${\cal Z}_n =\{1,2,\ldots ,n\equiv 0 \}$ of the residue classes by modulo $n$ .
In order to  have accordance with the marks of the elements of ${\cal S}_n$ and as $0\equiv n\; (\mbox{\rm mod}\; n)$, then the zero element in ${\cal Z}_n$ we denote $n$ (instead of 0). If $\pi \in {\cal S}_n$ and $\sigma\in {\cal S}_n$ are elements corresponding to the matrices $P\in {\cal P}_n$ and $S\in {\cal P}_n$ by the isomorphism of the groups ${\cal P}_n$ and ${\cal S}_n$ described by the formulas (\ref{pi}), (\ref{sigma}), (\ref{P}) and (\ref{S}), then for each $t=1,2,\ldots ,n$ there are:
\begin{equation}\label{l21}
(t)\pi = t\oplus 1 \quad (\mbox{by definition})
\end{equation}
\begin{equation}\label{l22}
(t)\sigma = n\oplus 1\ominus t\quad (\mbox{by definition})
\end{equation}

We will prove by induction, that for each whole positive number $l$ is true:
\begin{equation}\label{l23}
(t)\pi^l =t\oplus l
\end{equation}
When $l=1$, the proposition follows from (\ref{l21}). Let the equation $(t)\pi^l =t\oplus l$ be true. Then we get $(t)\pi^{l+1} =((t)\pi^l)\pi =(t\oplus l)\pi = t\oplus l\oplus 1$, and it follows that the equation (\ref{l23}) is true for every natural  number $l$.

Using the equations (\ref{l21}), (\ref{l22}) and (\ref{l23}) we continuously get:
$$(t)\pi^l \sigma =((t)\pi^l )\sigma =(t\oplus l)\sigma =n\oplus 1\ominus (t\oplus l)=n\oplus 1\ominus t \ominus l$$
$$(t)\sigma\pi^{n-l} = ((t)\sigma )\pi^{n-l} = (n\oplus 1\ominus t)\pi^{n-l} =(n\oplus 1\ominus t)\oplus (n-l) = $$
$$=2n\oplus 1\ominus t\ominus l = n\oplus 1\ominus t\ominus l$$

The last equation is true, because $2n\equiv n\equiv 0\; (\mbox{\rm mod}\; n)$. We see, that $(t)\pi^l \sigma =(t)\sigma\pi^{n-l}$ for every $t=1,2,\ldots ,n$ and therefore, $\pi^l \sigma =\sigma\pi^{n-l}$. Then having  in mind the isomorphism of the groups ${\cal S}_n$ and ${\cal P}_n$, it follows that the assertion of the lemma is true.

\hfill $\Box$

\begin{theor}\label{t1}
If $A\sim A^S$ and $B\sim A$ then $B\sim B^S$.
\end{theor}

Proof. Since $B\sim A$, then according to the lemma \ref{l1} there exist $k,l\in \{ 1,2,\ldots ,n\}$, such that $B=P^k AP^l$. Applying lemma \ref{l2} we get $BS=P^k AP^l S=P^k ASP^{n-l}$, and then follows, that $BS\sim AS$, i.e. according to (\ref{8}) $B^S \sim A^S$. But $A^S \sim A\sim B$ and because of the transitiveness of the relation $\sim$ we get $B^S \sim B\Rightarrow B\sim B^S$.

\hfill $\Box$

Theorem \ref{t1} gives us a reason to make the following definition:

\begin{dfn}\label{samoogledalna}
Let  $A\in {\cal Q}_n$. We say that $A$  is a representative of {\bf self-mirrored image} (or {\bf mirror image to itself})  interweaving, if
$A\sim A^S .$
\end{dfn}

Let us denote  the set $$\overline{{\cal M}_n} \subset \overline{{\cal Q}_n }$$ including all self-mirrored interweavings with repetition $n$.

For random $A\in {\cal B}_n$ we define the operation
\begin{equation}\label{A^R}
A^R =
\left(
\begin{array}{cccc}
a_{1\, n} & a_{2\, n} & \cdots & a_{n\, n} \\
a_{1\, n-1} & a_{2\, n-1} & \cdots & a_{n\, n-1} \\
\vdots    & \vdots    &        & \vdots    \\
a_{1\, 1} & a_{2\, 1} & \cdots & a_{n\, 1} \\
\end{array}
\right)
=\left( A^S \right)^T =\left( A S\right)^T =S^T  A^T = S A^T
\end{equation}

In other words the matrix $A^R$ is received by rotating the matrix $A$ by ${\rm 90}^\circ$ counter clockwise.

Obviously,
$$\left( \left( \left( A^R \right)^R \right)^R \right)^R =A .$$

In general  $A^R \ne A$.

\begin{lemma}\label{l3}
If $P$ is a binary matrix, defined by the formula (\ref{P}), then
$$P^T =P^{n-1} $$
\end{lemma}

Proof. If $P=(p_{ij} )$ and $P^T =(p_{i j}' )$, then by definition $p_{i j}' = p_{ji}$ for each $i,j\in \{ 1,2,\ldots n\}$. Let $PP^T =Q=(q_{ij} )\in {\cal P}_n$. Then for each $i=1,2,\ldots ,n$ there is $\displaystyle q_{ii} = \sum_{k=1}^n p_{ik} p_{ki}' =\sum_{k=1}^n p_{ik}^2 =(n-1)0 +1=1$ and it is the unique 1 in $i$th row of the matrix $Q=PP^T$. Therefore, $PP^T =E_n$, where $E_n$ is the identity matrix. We multiply from the left the two sides of the last equation with $P^{n-1}$ and have in mind, that $P^n =E_n$, then we get $P^{n-1} PP^T =P^{n-1} E_n$, and then finally we get, that $P^T =P^{n-1}$.

\hfill $\Box$

\begin{theor}\label{t2}
If $A\sim A^R$ and $B\sim A$ then $B\sim B^R$.
\end{theor}

Proof. $B\sim A$, hence according to the lemma \ref{l1} there exist natural numbers $k$ and $l$, as $B=P^k AP^l$. Then, when we apply lemma \ref{l2} and lemma \ref{l3} we get $B^R =SB^T =S(P^k A P^l )^T = S(P^T )^l A^T (P^T )^k =S(P^{n-1} )^l A^T (P^{n-1} )^k = SP^{nl-l} A^T P^{kn-k} = SP^{n-l}  A^T P^{n-k} =P^l S A^T P^{n-k}$. Therefore, $$B^R \sim A^R \sim A\sim B.$$
\hfill $\Box$

Theorem \ref{t2} gives us the right to give the following definition:
\begin{dfn}\label{rotat}
Let $A\in {\cal Q}_n$. If   $A \sim A^R$, then we say that $A$ is a representative of  {\bf rotation stable} interweaving.
\end{dfn}

Let denote   the set $$\overline{\cal R}_n \subset \overline{\cal Q}_n$$ of all rotation stable interweavings with repetition $n$.

The rotation stable interweavings play important role in practice.This means, that if a fabric is weaved which weaving structure is coded with a matrix, representative of rotation stable interweaving, then this fabric will have the same physical and mechanical characteristics (except of course the color) after a rotation by ${\rm 90}^\circ$.

\section{ A Representation of Binary Matrices}
 In \cite{umb2009} is described a representation of the elements of ${\cal B}_n$ using ordered $n$-tuples of natural numbers $\langle k_1 ,k_2 ,\ldots ,k_n \rangle$, where $0\le k_i \le 2^n -1$, $i=1,2,\ldots ,n$. The one to one corresponding is based on the binary presentation of the natural numbers.  It is not difficult to see that using this representation, this leads to making a fast and saving memory algorithms. Having in mind this we will create an algorithm, which will find just one representative of each equivalence class to the factor sets $\overline{{\cal Q}_n}$, $\overline{{\cal M}_n}$ and $\overline{{\cal R}_n}$. It will be the minimal element of the equivalence class with regard to the lexicographic order. As a consequence we will get an algorithm to solve the combinatoric problem to find the number of the equivalence classes in the sets ${\cal Q}_n$, ${\cal M}_n$ and ${\cal R}_n$ by the relation $\sim$ with given natural number $n$.

The matrices $P$ and $S$ given by the formulas (\ref{P}) and (\ref{S}) are coded using ordered $n$-tuples as it follows:
\begin{equation}\label{Ptuple}
P\; \equiv \; \langle2^{n-2},2^{n-3} ,\ldots ,2^1 ,2^0 ,2^{n-1} \rangle
\end{equation}
\begin{equation}\label{Stuple}
S\; \equiv \; \langle 2^0,2^1 , 2^2 ,\ldots  ,2^{n-2} ,2^{n-1} \rangle
\end{equation}

In programming languages  C, C++ and Java \cite{davis,Kernigan,Schildt} the number $x=2^k$ is calculated using a single  operation bitwise shift left ''$<<$'' with the statement $x=1<<k;$.

We examine the set ${\cal B} =\{ 0,1\}$. ${\cal B}$ together with the operations \textbf{conjunction} $\&\&$ (\textbf{AND}), \textbf{disjunction} $\|$ (\textbf{OR}) and \textbf{negation} $!$ (\textbf{NOT}) form the all known {\it boolean algebra} ${\cal B}(\&\& ,\| ,!)$, which role is very important in the computers and programming. We will use the pointed symbols for the operations in ${\cal B}(\&\& ,\| ,!)$, have in mind the semantic and the syntax of these operations in the widespread C++ programming language.

Let $A=(a_{ij} )$ and $B=(b_{ij} )$ are matrices of ${\cal B}_n$. Here and everywhere below, we will consider the semantic and syntax of C++ language and the first index will be 0, i.e $i,j\in \{ 0,1,2,\ldots ,n-1 \}$. We examine the following operations in ${\cal B}_n$, defined according to our aim as follows:

{\bf component conjunction}
\begin{equation}\label{1}
A \;\&\& \; B =C=(c_{ij} )
\end{equation}
where by definition for every $i,j \in \{ 0,1,2\ldots ,n-1 \} \qquad c_{ij} =a_{ij}  \;\&\& \; b_{ij}$

{\bf component disjunction}
\begin{equation}\label{2}
A \;\| \; B =C=(c_{ij} )
\end{equation}
where by definition for every $i,j \in \{ 0,1,2\ldots ,n-1 \}  \qquad c_{ij} =a_{ij}  \;\| \; b_{ij}$

{\bf component negation}
\begin{equation}\label{3}
! A =C=(c_{ij} )
\end{equation}
where by definition for every $i,j \in \{ 0,1,2\ldots ,n-1 \}  \qquad c_{ij} = ! a_{ij}$

{\bf transpose}
\begin{equation}\label{4}
t(A) =C=(c_{ij} )
\end{equation}
where by definition for every $i,j \in \{ 0,1,2\ldots ,n-1 \}  \qquad c_{ij} =a_{ji} $

{\bf logical product}
\begin{equation}\label{5}
A*B =C=(c_{ij} )
\end{equation}
where by definition for every $i,j \in \{ 0,1,2\ldots ,n-1 \} $
$$c_{ij} = \bigvee_{k=0}^{n-1} (a_{ik}  \; \&\& \; b_{kj} )=(a_{i\> 0}  \;\&\& \; b_{0\> j} )\; \|  \; (a_{i\> 1} \;\&\& \; b_{1\> j} ) \; \|  \; \cdots  \; \| \; (a_{i\> n-1}  \; \&\&  \; b_{n-1 \> j} )$$

It is not difficult to see that the introduced above binary operations are associative.

${\cal B}_n$ and the above-described operations $\&\&$, $\| $, $!$, $t()$ and $*$  form the algebra ${\cal B}_n (\&\& ,\| ,!,t(),*)$. Here and in the whole article the term {\it algebra} means {\it abstract algebra} considering the definition given in \cite{daintith}, and namely set equipped with various operations, assumed to satisfy some specified system of axiomatic laws. It is naturally we to put the linear order, and more precisely the lexicographic order in ${\cal B}_n (\&\& ,\| ,!,t(),*)$.

\section{The algebra ${\cal B}_n (\&\& ,\| ,!,t(),*) $ and bitwise operations.}

The use of bitwise operations is a well-working method, used in C/C++ and Java programming languages  \cite{davis,Kernigan,Schildt}.  Our aim is to show that the presentation of the binary matrices using ordered n-tuple of nonnegative integers and the bitwise operations make the realization of a better C++ class compared with the standard presentation of the binary matrices.

Let $x$, $y$ and $z$ are three integer variables, for which $w$ bits are needed. Let $x$ and $y$ are initialized and let $z=x\alpha x$, where $\alpha$ is one of the operators \textbf{bitwise conjunction (bitwise AND)}  \verb"&",  \textbf{bitwise disjunction (bitwise OR)} \verb"|"  or  \textbf{bitwise exclusive OR} \verb"^" . For each $i=0,1,\ldots ,n-1$ the new contents of the $i$ bit in $z$ will be as it is presented in  table \ref{bitwise}

bitwise conjunction $\&$, bitwise disjunction $|$, bitwise exclusive or \verb"^" and bitwise negation \verb"~"

\begin{table}[!h]
\begin{center}
\begin{tabular}{||c|c|c|c|c||}
  \hline\hline
The $i$ bit of  & The $i$ bit of  & The $i$ bit of  & The $i$ bit of & The $i$ bit of \\
 \verb"x" & \verb"y" & \verb"z = x & y;" & \verb"z = x | y;" & \verb"z = x ^ y;" \\
\hline\hline
  0 & 0 & 0 & 0 & 0 \\ \hline
  0 & 1 & 0 & 1 & 1 \\ \hline
  1 & 0 & 0 & 1 & 1 \\ \hline
  1 & 1 & 1 & 1 & 0 \\ \hline
  \hline
\end{tabular}
\caption{Bitwise operations in C++ programming languages}\label{bitwise}
\end{center}
\end{table}

In case \verb"z = ~x", if the $i$ bit of ix \verb"x" is 0, then the $i$ bit of \verb"z" becomes 1, and if the $i$ bit of \verb"x" is 1, then the $i$ bit of \verb"z" becomes 0, $i=0,1,\ldots ,w-1$.

If $k$ is a nonnegative integer, then the statement $z\; =\; x<<k$ (\textbf{bitwise shift left}) will write in the $(i+k)$-th  bit  of $z$ the value of the $k$-th  bit of $x$, where $i=0,1,\ldots ,w-k-1$, and the first $k$ bits of will be filled by zeroes. This operation is equivalent to a multiplication of $x$ by $2^k$. The statement $z\; =\; x>>k$  (\textbf{bitwise shift right}) work similarly.

For computer implementation of various algorithms related to the algebra ${\cal B}_n (\&\& ,\| ,!,t(),*) $, we propose to created the following C++ class:

\begin{verbatim}
class Bn_tuple {
    int n;
    int *Matr;
  public:
/* constructor without parameter: */
    Bn_X();
/* constructor with parameter n pointing the rows of the matrix: */
    Bn_X (unsigned int);
/*  copy constructor: */
    Bn_X (const Bn_X &);
/* destructor: */
    ~Bn_X();
/* predefines the operator "=": */
    Bn_X & operator = (const Bn_X &);
/* returns the size of the matrix: */
    int get_n() { return n; };
/* sets value 1 to the element (i,j) of the matrix: */
    void set_1 (int,int);
/* sets value 0 to the element (i,j) of the matrix: */
    void set_0 (int,int);
/* fills row i of the matrix using integer number r */
    void set_row(int,int);
/* gets element (i,j) of the matrix: */
    int get_element (int,int);
/* gets the row i of the matrix */
    int get_row (int);
/* predefines operators according to (1), (2), (3) and (5): */
    Bn_X operator && (Bn_X &);
    Bn_X operator || (Bn_X &);
    Bn_X operator ! ();
    Bn_X operator * (Bn_X &);
/* transposes matrix: */
    Bn_X t ();
/* defines order (lexicographical) */
    int operator < (Bn_X &);
}
\end{verbatim}

Furthermore, for our purposes, the class \verb"Bn_tuple" may be used for solving a number of tasks in the field of computer algebra\cite{tan}.

To evaluate the effectiveness and speed of the algorithms, which use objects of the algebra ${\cal B}_n (\&\& ,\| ,!,t(),*)$ it is necessary to evaluate  the algorithms, which realize the operations $\&\&$, $\|$ $!$, $*$, $t()$  and the operation ''$<$'' comparing two elements. In that sense we describe in details just these methods, realizing the above mentioned operations. To create these methods we use bitwise operations: bitwise conjunction $\&$, bitwise disjunction $|$, bitwise exclusive or \verb"^" and bitwise negation \verb"~", using these operations we raise the effectiveness and make the algorithms work faster. We suppose that the experienced programmer can easily create the other methods. Using the same model we can predefine the other relational operators: ''$<=$'', ''$>$'', ''$>=$'', ''$==$'' and ''$!=$''. If the dimensions of the two matrices, we compare, are not equal, then the relation ''$<$'' is not defined and the result we get is the negative number -1.

\begin{verbatim}
Bn_tuple Bn_tuple :: operator && (Bn_tuple &B) {
    Bn_tuple temp(n);
    if (B.get_n() != n)
      cout<<"unallowable value of a parameter \n";
    else
        for (int p=0; p<n; p++)
             *(temp.Matr + p) = *(this->Matr + p) & *(B.Matr + p);
    return temp;
}

Bn_tuple Bn_tuple :: operator || (Bn_tuple &B) {
    Bn_tuple temp(n);
    if (B.get_n() != n)
      cout<<"unallowable value of a parameter \n";
    else
        for (int p=0; p<n; p++)
             *(temp.Matr + p) = *(this->Matr + p) | *(B.Matr + p);
    return temp;
}

Bn_tuple Bn_tuple :: operator ! () {
	Bn_tuple temp(n);
	for (int i=0; i<n; i++) {
		 for (int j=0; j<n; j++) {
			 if ( get_element(i,j) ) temp.set_0(i,j);
				else temp.set_1(i,j);
		 }
	}
	return temp;
}

Bn_tuple Bn_tuple :: operator * (Bn_tuple &B) {
    Bn_tuple temp(n), TB(n);
    TB = t(B);
    int c, r_i, r_j;
    if (B.get_n() != n)
      cout<<"unallowable value of a parameter \n";
      else
      for (int i=0; i<n; i++)
        for (int j=0; j<n; j++) {
          r_i = this->get_row(i);
          r_j = TB.get_row(j);
          c = r_i & r_j;
          if (c==0) temp.set_0(i,j);
               else temp.set_1(i,j);
        }
    return temp;
}

Bn_tuple Bn_tuple :: t() {
	Bn_tuple temp(n);
	 int k;
	for (int i=0; i<n; i++)
		for (int j=0; j<n; j++) {
			k=get_element(i,j);
			if (k) temp.set_1(j,i);
				else temp.set_0(j,i);
		}
	return temp;
}

int Bn_tuple :: operator < (Bn_tuple &B) {
	int p = 0;
	if (B.get_n() != n)
	   cout<<"unallowable value\n";
	else
		while ((get_row(p)==B.get_row(p))&&(p<n-1 )) p++;
	if (get_row(p)<B.get_row(p) ) return 1;
			  else return 0;
}
\end{verbatim}

When we predefine operators $\&\&$, $||$ and $*$, if the dimensions of the two operands are not equal, then we receive the zero matrix of order the same as the first operand. But this result is not correct. We can say something more: in this case the operation is not defined, i.e. the result of the operation function is not correct and we have to be very careful in such situations. The situation is the same about the entered linear order, i.e. although the lexicographic order can be put for the words with different length, we examine by definition only the matrices of one and the same order. The result we receive when we compare the matrices with different dimensions is not correct. We give suitable massages in this situations.

It is easy to convince that the following proposition is true:

\begin{prop} \label{p2}
For each nonnegative  integer $n$, for the computer representation via the class \verb"Bn_tuple", using the C++ programming language, of the algebra ${\cal B}_n (\&\& ,\| ,!,t(),*)$, the following assertions are true:

(i) For each of the operations $\&\&$, $\|$, $<$ are needed in $O(n)$ standard operations;

(ii) The transpose and negation operation require $O(n^2)$ standard operations, each;

(iii) When we use the above-described realization of the operation functions $*$, this operation requires $O(n^2)$ standard operations;

(iv) For every object of the class \verb"Bn_tuple" are necessary $O(n)* {\rm sizeof}\> {\rm(int)}$ bytes of the operating memory;

(v) Initialization of an object of the class \verb"Bn_tuple" requires $O(n)$ standard operations.

\hfill $\Box$
\end {prop}

If the binary matrix $A$ is represented using the ordered $n$-tuple of integers, then to check whether $A$ belongs to the set ${\cal Q}_n \subset {\cal B}_n$ we can use the following obvious assertion:

\begin{lemma}\label{l4}
 Let $A\in {\cal B}_n$ be represented by the ordered $n$-tuple $\langle k_1 ,k_2 ,\ldots ,k_n \rangle$, where $0\le k_i \le 2^n -1$, $i=1,2,\ldots ,n$. Let $|$ and $\&$ denote the operations bitwise OR and bitwise AND respectively. Then:

(i) The number $k_i$, $i=1,2,\ldots ,n$ represents row of zeroes if and only if $k_i =0$;

(ii) The number $k_i$, $i=1,2,\ldots ,n$ represents row of ones if and only if $k_i =2^n -1$;

(iii) $j$th column of $A$ is made of zeroes if and only if $$(k_1 \, | \, k_2 \, | \cdots |\, k_n )\, \& \, 2^j =0;$$

(iv) $j$th column of $A$ is made of ones if and only if $$(k_1 \, \& \, k_2 \, \& \cdots \& \, k_n )\, \& \, 2^j \ne 0.$$
\hfill $\Box$
\end{lemma}

The algorithm, which is below-described is based on the following assertions:
\begin{lemma}\label{l5}
If $A\in {\cal B}_n $, $B\in {\cal P}_n$, then
$$A*B=AB$$ and
$$B*A=BA$$
\end{lemma}

Proof. Let $A=(a_{ij} )$, $B=(b_{ij} )$, $U=A*B=(u_{ij} )$ and $V=AB=(v_{ij} )$, $i,j=1,2,\ldots ,n$. Let the unique 1 in the $j$th column of $B\in {\cal P}_n$ be on the $s$th place, i.e. $b_{sj} =1$ and $b_{kj} =0$ when $k\ne s$. Then by definition
$$u_{ij} =\bigvee_{k=1}^n (a_{ik}  \; \& \; b_{kj} )=
\left\{
\begin{array}{ccc}
1 & \mbox{\rm for}  & a_{is} =1\\
0 & \mbox{\rm for}  & a_{is} =0
\end{array}
\right.$$ and
$$v_{ij} =\sum_{k=1}^n (a_{ik}  b_{kj} )=
\left\{
\begin{array}{ccc}
1 & \mbox{\rm for}  & a_{is} =1\\
0 & \mbox{\rm for}  & a_{is} =0
\end{array}
\right.$$

Therefore, $u_{ij} = v_{ij}$ for each $i,j\in \{ 1,2,\ldots , n\}$.

Analogously can be proved, that $B*A=BA$.

\hfill $\Box$

\begin{lemma}\label{l6}
 Let $A\in {\cal B}_n$ be represented using ordered $n$-tuple $\langle k_1 ,k_2 ,\ldots ,k_n \rangle$ and let $A$ be a minimal element of the equivalence class corresponding to the lexicographic order in $\mathbb{N}^n$. Then $k_1 \le k_t$ for each $t=2,3,\ldots ,n$.
\end{lemma}

Proof. We presume, that there exists $t\in \{ 2,3,\ldots ,n\}$, such as $k_t <k_1$. Then if we move the first row on the last place $t-1$ times, we get a matrix $A' \in{\cal B}_n$, such as $A' \sim A$ and $A'$ is represented using the $n$-tuple $\langle k_t ,k_{t+1} ,\ldots ,k_n ,k_1 ,\ldots ,k_{t-1} \rangle$. Then obvious $A' <A$ according to the lexicographic order in $\mathbb{N}^n$, which runs counter to the minimum of $A$ in the equivalence class $\overline{A}$.

\hfill $\Box$

We can create the following generalized algorithm to obtain a representative of each equivalence class in the factor sets $\overline{\cal Q}_n$, $\overline{\cal M}_n$ and  $\overline{\cal R}_n$

\begin{algr}\label{alg1}.
\begin{enumerate}
\item\label{it1} Generate all ordered $n$-tuples of natural numbers $\langle k_1 ,k_2 ,\ldots ,k_n \rangle$ such as $1\le k_i \le 2^n -2 $ and $k_1 \le k_i$ as $i=2,3,\ldots ,n$;
\item\label{it2} Check if the elements obtained in \ref{it1} belong to the set ${\cal Q}_n$ according to the lemma \ref{l4} (iii)  (iv) (cases (i) and (ii) we reject when we generate the elements in point \ref{it1} according to the lemma \ref{l6});
\item\label{it3} Check whether the element, obtained in point \ref{it2} is minimal in the equivalence class. According to the lemmas \ref{l1} and \ref{l5} $A$ is minimal in $\overline{A}$ if and only if $A\le P^k *A*P^l$ for each $k,l\in \{ 0,1,\ldots ,n-1\}$, where $P$ is the matrix represented by $n$-tuple (\ref{Ptuple});
\item\label{it4} Check whether the elements obtained in point \ref{it3} belong to the set ${\cal M}_n$ according to definition \ref{samoogledalna} and applying lemmas \ref{l1} and \ref{l5};
\item\label{it5} Check whether the elements obtained in point \ref{it3} belong to the set ${\cal R}_n$ according to definition  \ref{rotat} and applying lemmas \ref{l1} and \ref{l5}.
\end{enumerate}
\hfill $\Box$
\end{algr}

\section{Conclusion}

In this work  we proved that there are algorithms which need $O(n)$ standard operations and implement the operations $\&\&$, $\|$, $<$ and
there is an algorithm which needs $O(n^2 )$ operations and implements the operation logical multiplication $*$ of two binary matrices, which are represented using ordered $n$-tuple of integers. In the same time, as is well known, to implement the operations $\&\&$, $\|$, $<$ according to the classical definition we need $O(n^2 )$ standard operations and to implement the operation logical multiplication $*$ of two integer  matrices  according to the classical definition we need $O(n^3 )$ standard operations. For every object of the described in this work C++ class \verb"Bn_tuple" are necessary $O(n)* {\rm sizeof}\> {\rm(int)}$ bytes of the operating memory. For every integer matrix are necessary $O(n^2 )* {\rm sizeof}\> {\rm(int)}$ bytes of the operating memory. Finally we saw that initialization of an object of the class \verb"Bn_tuple" requires $O(n)$ standard operations whereas to initialize a $n\times n$ integer matrix $O(n^2 )$ standard operations are necessary. This proved the usefulness of the use of bitwise operations in programming.

Applying the above ideas, a computer program that receives some different weaving structures is created.
The work results of this program taking some values of $n$ are generalized in the table \ref{tabl2}

\begin{table}[!h]
\begin{center}
\begin{tabular}{||c||c|c|c|c|c||}
\hline\hline
$n$                                      &   2  &   3   &   4       &      5   & 6   \\
\hline\hline
$\displaystyle |{\cal B}_n |$   &   16  &   512  &  65 536    & 33 554 432  &  $2^{36} > 2^{32} -1$   \\
\hline
$\displaystyle |{\cal Q}_n |$   &   2  &   102  &  22 874    & 17 633 670  &   $>2^{32} -1$   \\
\hline
$\displaystyle |\overline{\cal B}_n |$   &   7  &   64  &  4 156    & 1 342 208  & 1 908 897 152     \\
\hline
$\displaystyle |\overline{\cal Q}_n |$   &   1  &   14  &  1 446    & 705 366  &  1 304 451 482    \\
\hline
$\displaystyle |\overline{\cal M}_n |$   &   1  &   2   &    142    &   1 302  &   586 060    \\
\hline
$\displaystyle |\overline{\cal R}_n |$   &   1  &   2   &     18    &      74  &    902   \\
\hline\hline
\end{tabular}
\caption{}\label{tabl2}
\end{center}
\end{table}

When $n\ge 6$ there are got too large values (see formula (\ref{bn})) and to avoid ''overloading'' it is necessary to be used some special programming techniques\cite{tan} which is not the task in this work.


\begin {thebibliography}{99}


\bibitem{bogopolski}
\textsc{O. Bogopolski} Introduction to Group Theory",
\emph{European Mathematical Society, Zurich} 2008.

\bibitem{borzunov} \textsc{G. I. Borzunov } Textile Industry - Survey Information.
\emph{Moscow, CNII ITEILP,} vol. 3, 1983 (in Russian).

\bibitem{daintith} \textsc{J. Daintith, R. D. Nelson} The Penguin Dictionary of Mathemathics.
\textit{Penguin books}, 1989.

\bibitem{davis} \textsc{S. R. Davis} C++ for Dummies.
\emph{IDG Books Worldwide,} 2000.

\bibitem{Kernigan} \textsc{B. W. Kernigan, D. M Ritchie}  The C Programming Language.
\emph{ AT$\&$T Bell Laboratories,}  1998.

\bibitem{Schildt} \textsc{H. Schildt} Java 2 A Beginner's Guide.
\emph{ McGraw-Hill}, 2001.

\bibitem{sachkov} \textsc{V. N. Sachkov, V. E. Tarakanov}  Combinatorics of Nonnegative Matrices.
\emph{Amer. Math. Soc.,}  1975.

\bibitem{serre} \textsc{J.-P. Serre} Linear Representations of Finite Groups.
\emph{Springer-Verlag, New York}, 1977.

\bibitem{tan} \textsc{Tan Kiat Shi, W.-H. Steeb, Y. Hardy} Symbolic C++: An Introduction to Computer Algebra using Object-Oriented Programming.
\textit{Springer}, 2001.

\bibitem{tarakanov} \textsc{V. E. Tarakanov} Combinatorial Problems and (0,1)-matrices.
\emph{Moscow, Nauka,}  1985 (in Russian).

\bibitem{vos} \textsc{A. De Vos} Reversible Computing: Fundamentals, Quantum Computing, and Applications.
\emph{Wiley}, 2010.

\bibitem{yordzhev} \textsc{K. Yordzhev, I. Statulov}  Mathematical Modeling and Quantitative Evaluation of Primary Weaving Braids.
\emph{Textiles and Clothing},  10, (1999), 18--20 (in Bulgarian).

\bibitem{bansko} \textsc{K. Yordzhev} On an Equivalence Relation in the Set of the Permutation Matrices.
\emph{Blagoevgrad, Bulgaria, SWU, Discrete Mathematics and Applications},  (2004), 77--87.

\bibitem{umb2009} \textsc{K. Yordzhev} An Example for the Use of Bitwise Operations in programming.
\emph{Mathematics and education in mathematics},  38 (2009), 196--202.

\bibitem{KYKo} \textsc{K. Yordzhev, H. Kostadinova} Using Mathematical Methods in the Weaving for Receiving Quantity Valuations of the Textile Structure Variety
\emph{Textiles and Clothing},  1, (2011), 7--10  (in Bulgarian).

\end{thebibliography}
\end{document}